\documentclass{elsart}
\usepackage{amssymb}
\usepackage{graphicx}

\begin{document}

\begin{frontmatter}

\title{Global bifurcations of limit cycles
in~the~Kukles~cubic~system\thanksref{label1}}
\thanks[label1]{This work was supported by grants from the German Academic Exchange
Service (DAAD) and the Netherlands Organization for Scientific Research (NWO).
The author is very grateful to the Department of Ma\-the\-ma\-tics,
Technische Universit\"{a}t Dresden (Germany) for hospitality during his
stay in September\,--\,November 2016 and thanks Prof. Stefan Siegmund
for kindly support and challenging discussions.}

\author{Valery A. Gaiko}

\ead{valery.gaiko@gmail.com}

\address{National~Academy~of~Sciences~of~Belarus,
United~Institute~of~Informatics~Problems,
Minsk~220012,~Belarus}

\begin{abstract}
In this paper, using our bifurcational geometric approach, we solve the problem
on the maximum number and distribution of limit cycles in the Kukles system
representing a planar polynomial dynamical system with arbitrary linear and
cubic right-hand sides and having an anti-saddle at the origin. We also apply
alternatively the Wintner--Perko termination principle to solve this problem.
\end{abstract}

\end{frontmatter}

\section{Introduction}

In this paper, we continue studying the Kukles cubic system
	$$
	\dot{x}=y,\quad
	\dot{y}=-x+\delta y+a_1x^2+a_2xy+a_3y^2+a_4x^3+a_5x^2y+a_6xy^2+a_7y^3\!.~~
	\eqno(1.1)
	$$
I.\,S.\,Kukles was the first who began to study~(1.1) solving
the center-focus problem for this system in 1944: he gave
the necessary and sufficient conditions for $O(0,0)$ to be a center for~(1.1)
with $a_7=0$~\cite{isk}. Later, system~(1.1) was studied by many mathematicians.
For example, in~\cite{lp} the necessary and sufficient center conditions for arbitrary
system~(1.1), when $a_7\neq0,$ were conjectured. In~\cite{rst}, global qualitative pictures
and bifurcation diagrams of a reduced Kukles system $(a_7=0)$ with a center were given.
In~\cite{wcy}, the global analysis of system~(1.1) with two weak foci was carried out.
In~\cite{yeye}, the number of singular points under the conditions of a center or a weak
focus for~(1.1) was investigated. In \cite{zzhtt}, new distributions of limit cycle for
the Kukles system were obtained. In \cite{rr}, an accurate bound of the maximum number
of limit cycles in a class of Kukles type systems was provided.
  \par
In \cite{gvh04,g08b}, we constructed a canonical cubic dynamical system of Kukles type
and carried out the global qualitative analysis of a special case of the Kukles system
corresponding to a generalized cubic Li\'{e}nard equation. In particular, it was shown that
the foci of such a Li\'{e}nard system could be at most of second order and that such system
could have at most three limit cycles in the whole phase plane. Moreover, unlike all previous
works on the Kukles type systems, global bifurcations of limit and separatrix cycles using
arbitrary (including as large as possible) field rotation parameters of the canonical system
were studied. As a result, a classification of all possible types of separatrix cycles for
the generalized cubic Li\'{e}nard system was obtained and all possible distributions of
its limit cycles were found.
  \par
In \cite{G,g05,g08a,g09b}, we also presented a solution of Hilbert's sixteenth problem in
the quadratic case of polynomial systems proving that for quadratic systems four is really
the maximum number of limit cycles and $(3\!:\!1)$ is their only possible distribution.
We established some preliminary results on gene\-ralizing our ideas and methods
to special cubic, quartic and other polynomial dynamical systems as~well.
In~\cite{g08b,g09a,g11b,g12a}, e.\,g., we presented a solution of Smale's thirteenth
problem \cite{sml} proving that the classical Li\'{e}nard system with a polynomial of
degree $2k+1$ could have at most $k$ limit cycles and we could conclude that
our results agree with the conjecture of~\cite{ldmp} on the maximum number
of limit cycles for the classical Li\'{e}nard polynomial system. In~\cite{g12b,g12c,g14},
under some assumptions on the parameters, we found the maximum number of limit cycles
and their possible distribution for the general Li\'{e}nard polynomial system. In~\cite{g11a},
we studied multiple limit cycle bifurcations in the well-known FitzHugh--Nagumo neuronal model.
In~\cite{bg,g16}, we completed the global qualitative analysis of quartic dynamical systems
which model the dynamics of the populations of predators and their prey in a given ecological
system.
  \par
In this paper, we will use the obtained results and our bifurcational geometric approach for
studying limit cycle bifurcations of Kukles cubic system (1.1). In~Section~2, we construct
new canonical systems with field rotation parameters for studying global bifurcations of
limit cycles of (1.1). In~Section~3, using these canonical systems and geometric properties
of the spirals filling the interior and exterior domains of limit cycles, we give a solution
of the problem on the maximum number and distribution of limit cycles for Kukles system (1.1).
In~Section~4, applying the Wintner--Perko termination principle, we give an alternative solution
of this problem. This is related to the solution of Hil\-bert's sixteenth problem on the maximum
number and distribution of limit cycles in planar polynomial dynamical systems \cite{G}.

\section{Canonical Systems}

Applying Erugin's two-isocline method \cite{G} and studying the rotation properties \cite{BL,G,P}
of all the parameters of~$(1.1),$ we prove the following theorem.
    \par
    \medskip
\noindent \textbf{Theorem 2.1.}
\emph{Kukles system $(1.1)$ with limit cycles can be reduced to the canonical form}
    \vspace{-4mm}
	  $$
    \begin{array}{c}
    \dot{x}=y
    \equiv P(x,y),\\[1mm]
    \dot{y}=q(x)\!+\!(\alpha_0-\beta+\gamma+\beta\,x+\alpha_2\,x^2)\,y\!
    +\!(c+dx)\,y^2+\gamma\,y^3
    \equiv Q(x,y),
    \end{array}	
	  \eqno(2.1)
	  $$
\emph{where}

$1)$~$q(x)=-x+(1+1/a)\,x^2-(1/a)\,x^3\!,\; a=\pm1,\pm2\;$ \emph{or}

$2)$~$q(x)=-x+b\,x^3\!,\; b=0,-1,$ \emph{or}

$3)$~$q(x)=-x+x^2;$

\noindent \emph{$\alpha_0,$ $\alpha_2,$ $\gamma$ are field rotation parameters and
$\beta$ is a~semi-rotation parameter.}
	\medskip
	\par
\noindent \textbf{Proof.} \ System (1.1) has two basic isoclines: the cubic curve
	$$
	-x+\delta y+a_1x^2+a_2xy+a_3y^2+a_4x^3+a_5x^2y+a_6xy^2+a_7y^3=0
	$$
as the isocline of ``zero'' and the straight line $y=0$ as the isocline of ``infinity''.
	\par
These isoclines intersect at least at one point: at the origin which is an anti-saddle
(a center, a focus or a node). Besides, (1.1) can have two more finite singularities
(two saddles or a saddle and an anti-saddle) or one additional finite singular point
(a saddle or a saddle-node), or no other finite singularities at all. All these singular
points lie on the $x$-axis $(y=0),$ and their coordinates are defined by the equation
	$$
	q(x)\equiv-x+a_1x^2+a_4x^3=0
	\eqno(2.2)
	$$
\noindent depending just on the parameters $a_1$ and $a_4.$
	\par
Without loss of generality, $q(x)$ as given by (2.2) can be written in the following forms:

$1)$~$q(x)\equiv-(1/a)x(x\!-\!1)(x\!-\!a)=-x\!+\!(1\!+\!1/a)\,x^2\!-\!(1/a)\,x^3\!,\; a=\pm1,\pm2\;$ or

$2)$~$q(x)\equiv-x(1-bx^2)=-x+b\,x^3\!,\; b=0,-1,$ or

$3)$~$q(x)\equiv-x(1-x)=-x+x^2.$
	\par
It means that, together with the anti-saddle in $(0,0),$ we can have also:

$1)$~two saddles: at $(1,0)$ and $(-2,0)$ for $a=-2$ or at $(1,0)$ and $(-1,0)$ for $a=-1;$
or a~saddle at $(1,0)$ and an anti-saddle at $(2,0)$ for $a=2;$ or a~saddle-node at $(1,0)$
for $a=1;$

$2)$~no other finite singularities;

$3)$~one saddle at $(1,0).$
	\par
At infinity, system (1.1) has at most four singular points: one of them is in the vertical
direction and the others are defined by the equation
	$$
	a_7u^3+a_6u^2+a_5u+a_4=0,\quad u=y/x.
	\eqno(2.3)
	$$
Instead of the parameters $\delta,$ $a_2,$ $a_3,$ $a_5,$ $a_6,$ $a_7,$
also without loss of generality, we can introduce some new parameters $c,$~$d,$~$\alpha_0,$~$\alpha_2,$~$\beta,$~$\gamma\!:$
	$$
	\delta=\alpha_0-\beta+\gamma;\;\; a_2=\beta;\;\; a_3=c;\;\;
	a_5=\alpha_2;\;\; a_6=d;\;\; a_7=\gamma
	$$
\noindent to have more regular rotation of the vector field of~(1.1)~\cite{G}.
	\par
Then, taking into account $q(x),$ equation (2.3) is written in the form
	$$
	\gamma\,u^3+d\,u^2+\alpha_2\,u+s=0,\quad u=y/x,\quad s=-1/a,b.
	\eqno(2.4)
	$$
Thus, we have reduced (1.1) to canonical system (2.1).
	\par
If $c=d=\alpha_0=\alpha_2=\beta=\gamma=0,$ we obtain the following Hamiltonian systems:
    \par
\vspace{-4mm}
	$$
	\dot{x}=y,\quad
	\dot{y}=-x+(1+1/a)\,x^2-(1/a)\,x^3\!,\quad a=\pm1,\pm2;
	\eqno(2.5)
	$$
\vspace{-2mm}
	$$
	\dot{x}=y,\quad
	\dot{y}=-x+b\,x^3\!,\quad b=0,-1;
	\eqno(2.6)
	$$
\vspace{-2mm}
	$$
	\dot{x}=y,\quad
	\dot{y}=-x+x^2.
	\eqno(2.7)
	$$
    \par
All their vector fields are symmetric with respect to the $x$-axis, and, besides, the fields
of system~(2.5) with $a=2,-1$ and system~(2.6) with $b=0,-1$ are symmetric with respect
to the straight line $x=1$ and centrally symmetric with respect to the point $(1,0).$
Systems (2.5)--(2.7) have the following Hamiltonians, respectively:
    \par
\vspace{-4mm}
	$$
	H(x,y)=x^2-(2/3)\,(1+1/a)\,x^3+(1/(2a))\,x^4+y^2\!,\quad a=\pm1,\pm2;
	$$
\vspace{-2mm}
	$$
	H(x,y)=x^2-(b/2)\,x^4+y^2\!,\quad b=0,-1;
	$$
\vspace{-2mm}
	$$
	H(x,y)=x^2-(2/3)\,x^3+y^2\!.
	$$
    \par
If $\alpha_0=\alpha_2=\beta=\gamma=0,$ we will have the system
	$$
	\dot{x}=y,\quad
	\dot{y}=q(x)+(c+dx)\,y^2
	\eqno(2.8)
	$$
and the corresponding equation
  $$
  \frac{dy}{dx}
  =\frac{q(x)+(c+dx)\,y^2}{y}
  \equiv F(x,y).
	\eqno(2.9)
	$$
Since $F(x,-y)=-F(x,y),$ the direction field of (2.9) (and the vector field of (2.8)
as~well) is symmetric with respect to the $x$-axis. It follows that system~(2.8) has
centers as anti-saddles and cannot have limit cycles surrounding these points.
Therefore, without loss of generality, the parameters $c$ and $d$ in system (2.1)
can be fixed.
    \par
To prove that the parameters $\alpha_{0},$ $\alpha_{2},$ $\gamma$ and $\beta$ 
rotate the vector field of~(2.1), let~us calculate the following determinants:
    $$
    \begin{array}{c}
    \Delta_{\alpha_0}=P\,Q'_{\alpha_0}-QP\,'_{\alpha_0}=y^{2}\geq0,\\[1mm]
    \Delta_{\alpha_2}=P\,Q'_{\alpha_2}-QP\,'_{\alpha_2}=x^{2}y^{2}\geq0,\\[1mm]
    \Delta_{\gamma}=P\,Q'_{\gamma}-QP\,'_{\gamma}=y^{2}(1+y^{2})\geq0,\\[1mm]
    \Delta_{\beta}=P\,Q'_{\beta}-QP\,'_{\beta}=(x-1)\,y^{2}.
    \end{array}
    \\[-4mm]
    $$
    \par
By definition of a field rotation parameter~\cite{BL,G}, for increasing
each of the parameters $\alpha_0,$ $\alpha_2$ and $\gamma,$
under the fixed others, the vector field of system~(2.1) is rotated in
positive direction (counter\-clock\-wise) in the whole phase plane;
and, conversely, for decreasing each of these parameters, the vector
field of~(2.1) is rotated in negative direction (clock\-wise).
For increasing the parameter $\beta,$ under the fixed others, the vector field of
system~(2.1) is rotated in~positive direction (counter\-clock\-wise) in the
half-plane $x>1$ and in negative direction (clock\-wise) in~the half-plane
$x<1$ and vice versa for decreasing this parameter.
We will call such a parameter as a semi-rotation one.
    \par
Thus, for studying limit cycle bifurcations of~(1.1), it is sufficient to consider
canonical system~(2.1) containing the field rotation para\-me\-ters
$\alpha_0,$ $\alpha_2,$ $\gamma$ and the semi-rotation parameter
$\beta.$ The theorem is proved.\quad$\Box$

\section{Global Bifurcations of Limit Cycles}

By means of our bifurcational geometric approach
\cite{gvh04,g08b,g09a,g11b,g12a,g12b,g12c,g14,g15,g16},
we will consider now the Kukles cubic system in the form (when $a=2)\!:$
	  $$
    \begin{array}{c}
    \dot{x}=y,\\[1mm]
    \dot{y}=-(1/2)x(x-1)(x-2)\!+\!(\alpha_0-\beta+\gamma+\beta\,x+\alpha_2\,x^2)\,y\!\\[1mm]
    +\,(c+dx)\,y^2+\gamma\,y^3\!.
    \end{array}	
	  \eqno(3.1)
	  $$
    \par
All other Kukles systems can be considered in a similar way. Using system (3.1),
we will prove the following theorem.
    \par
\medskip
\noindent \textbf{Theorem 3.1.}
\emph{Kukles cubic system $(1.1)$ can have at most four limit cycles in~$(3\!:\!1)$-distribution.}
    \par
\medskip
\noindent \textbf{Proof.} \
According to Theorem~2.1, for the study of limit cycle bifurcations of system~(1.1),
it is sufficient to consider canonical system~(2.1) containing the field rotation
para\-me\-ters $\alpha_0,$ $\alpha_2,$ $\gamma$ and the semi-rotation parameter
$\beta.$ We~will work with system (3.1) which has three finite singularities:
a saddle $S(1,0)$  and two anti-saddles, $O(0,0)$ and $A(2,0).$
    \par
Vanishing all of the rotation parameters $\alpha_0,$ $\alpha_2,$ $\gamma$ and
also the parameter $\beta,$ we will get the system
    $$
    \begin{array}{c}
    \dot{x}=y,\quad
    \dot{y}=-(1/2)x(x-1)(x-2)+(c+dx)\,y^2
    \end{array}
	  \eqno(3.2)
	  $$
which is symmetric with respect to the $x$-axis and has centers as anti-saddles
at the points $O(0,0)$ and $A(2,0).$ Its center domains are bounded by separatrix
loops of the saddle $S(1,0).$
    \par
Let us input successively the field rotation parameters into (3.2).
Begin with the parameter~$\alpha_0$ supposing that $\alpha_0>0\!:$
    $$
    \begin{array}{c}
    \dot{x}=y,\quad
    \dot{y}=-(1/2)x(x-1)(x-2)+\alpha_0\,y+(c+dx)\,y^2.
    \end{array}
    \eqno(3.3)
    $$
On increasing $\alpha_0,$ the vector field of~(3.3) is rotated in positive
direction (counter\-clockwise) and the centers $O$ and $A$ turn into unstable foci.
    \par
Fix $\alpha_0$ and input the parameter $\beta>0$ into (3.3):
    $$
    \begin{array}{c}
    \dot{x}=y,\quad
    \dot{y}=-(1/2)x(x-1)(x-2)+(\alpha_0-\beta+\beta x)\,y+(c+dx)\,y^2.
    \end{array}
    \eqno(3.4)
    $$
Then, in the half-plane $x>1,$ the vector field of~(3.4) is rotated
in positive direction again and the focus $A$ remains unstable;
in the half-plane $x<1,$ the vector field is rotated in negative
direction and, when $\beta=\alpha_0>0,$ the focus $O$ becomes weak.
Fix this value of the parameter $\beta=\beta^{AH}$ (the Andronov--Hopf
bifurcation value).
    \par
Fix the parameters $\alpha_0>0,$ $\beta=\beta^{AH}>0$ and input
the third parameter, $\alpha_2<0,$ into this system:
    \vspace{-4mm}
    $$
    \begin{array}{c}
    \dot{x}=y,\\[1mm]
    \dot{y}=-(1/2)x(x-1)(x-2)+(\alpha_0-\beta+\beta x+\alpha_2 x^2)\,y+(c+dx)\,y^2.
    \end{array}
    \eqno(3.5)
    $$
The vector field of~(3.5) is rotated in negative direction
(clockwise) and a~big stable limit cycle appears immediately
from infinity. Denote this cycle by~$\Gamma_1^{bc}.$
    \par
On decreasing $\alpha_2,$ the cycle $\Gamma_1^{bc}$ will contract and,
for some value $\alpha_2=\alpha_2^{8l},$ a~separatrix eight-loop of the saddle $S$
will be formed around the points $O$~and~$A.$ On further decreasing $\alpha_2,$
two stable limit cycle, $\Gamma_1^{O}$ and $\Gamma_1^{A},$ will appear from the
eight-loop surrounding $O$~and~$A,$ respectively. These cycles will
contract and, finally, will disappear at the foci $O$~and~$A.$
    \par
Suppose that on decreasing $\alpha_2,$ the limit cycle $\Gamma_1^O$
and $\Gamma_1^A$ still exist and consider logical possibilities
of the appearance of other (semi-stable) limit cycles from a
``trajectory concentration'' surrounding the points $O$~and~$A.$
    \par
Denote the domains outside the cycle $\Gamma_1^O$ and $\Gamma_1^A$
by $D_1^O$ and $D_1^A$, the domains inside the cycles by $D_2^O$
and $D_2^A,$ respectively. It is clear that on decreasing
$\alpha_2,$ a semi-stable limit cycle cannot appear in the domains
$D_1^O$ and $D_1^A,$ since the focus spirals filling these domains will
untwist and the distance between their coils will increase because
of the vector field rotation in negative direction.
    \par
By contradiction, we can also prove that a semi-stable limit cycle
cannot appear in the domains $D_2^O$ and $D_2^A.$ Suppose it appears in
a domain for some values of the parameters: $\alpha_0^*>0,$ $\alpha_2^*<0,$
$\beta^{AH}>0.$ Return to initial system (3.2) and change the
order of inputting the field rotation parameters. 
    \par
Input first the parameter $\alpha_2<0\!:$
    $$
    \begin{array}{c}
    \dot{x}=y,\quad
    \dot{y}=-(1/2)x(x-1)(x-2)+\alpha_2\,x^2y+(c+dx)\,y^2.
    \end{array}
    \eqno(3.6)
    $$
Fix it under $\alpha_2=\alpha_2^*.$ The vector field of~(3.6) is rotated
in negative direction and the points $O$ and $A$ become stable foci.
    \par
Inputting the parameter $\beta>0$ into (3.6), we will have the system
    \vspace{-2mm}
    $$
    \begin{array}{c}
    \dot{x}=y,\\[1mm]
    \dot{y}=-(1/2)x(x-1)(x-2)+(-\beta+\beta x+\alpha_2\,x^2)\,y+(c+dx)\,y^2,
    \end{array}
    \eqno(3.7)
    $$
the vector field of which is rotated in positive direction in the
half-plane $x>1$ and in negative direction in the half-plane $x<1.$
Fix it under $\beta=\beta^{AH}\!.$
    \par
Inputting the parameter $\alpha_0>0$ into (3.7), we will get again
system (3.5), where the vector field is rotated in positive direction.
Under this rotation, stable limit cycles, $\Gamma_1^O$ and $\Gamma_1^A,$
will appear from the foci $O$ and $A,$ when they change the character of stability.
These cycles will expand, the focus spirals will untwist and the distance between
their coils will increase on increasing the parameter $\alpha_0$ to the
value $\alpha_0=\alpha_0^*.$ It follows that there are no values of
$\alpha_0=\alpha_0^*>0,$ $\alpha_2=\alpha_2^*<0$ and $\beta=\beta^{AH}>0,$ 
for which a semi-stable limit cycle could appear in the domains $D_2^O$
and $D_2^A.$
    \par
Thus, we have proved the uniqueness of limit cycles surrounding the points
$O$ and $A$ for $\alpha_0>0,$ $\alpha_2<0$ and $\beta=\beta^{AH}>0.$
Similarly, it can be proved the uniqueness of a big limit cycle
surrounding all the finite singularities $O,$ $A$ and $S$
for this set of the parameters.
    \par
Consider again system (3.5) for $\alpha_0>0,$ $\alpha_2<0$ and $\beta=\beta^{AH}>0$
supposing that it has two stable limit cycles, $\Gamma_1^O$ and $\Gamma_1^A.$
Change the parameter $\beta\!:$ $\beta>\beta^{AH}=\alpha_0>0.$ On increasing this
parameter, the weak focus $O$ will become rough stable one generating an unstable
limit cycle, $\Gamma_2^O$ (the Andronov--Hopf bifurcation). On further increasing
$\beta,$ the limit cycle $\Gamma_2^O$ will join with $\Gamma_1^O$ forming a semi-stable
limit cycle, $\Gamma_{12}^O,$ which will disappear in a ``trajectory concentration''
surrounding the point $O.$ Can another semi-stable limit cycle appear around this point
in addition to $\Gamma_{12}^O?$ It is clear that such a limit cycle cannot appear either
in the domain $D_3^O$ bounded by the origin $O$ and $\Gamma_2^O$ or in the domain $D_1^O$
bounded on the inside by $\Gamma_1^O$ because of the increasing distance between the
spiral coils filling these domains under increasing $\beta.$
    \par
To prove impossibility of the appearance of a semi-stable limit
cycle in the domain $D_2^O$ bounded by the cycles $\Gamma_1^O$
and $\Gamma_2^O$ (before their joining), suppose the contrary,
i.\,e., for some set of values of the parameters $\alpha_0^*>0,$ 
$\alpha_2^*<0$ and $\beta^*>0,$ such a semi-stable cycle exists.
Return to system (3.2) again and input the parameters
$\alpha_2<0$ and $\beta>0$ getting system (3.7).
In the half-plane $x<1,$ both parameters act in a similar way:
they rotate the vector field of (3.7) in negative direction
turning the origin $O$ into a stable focus. In the half-plane
$x>1,$ they rotate the field in opposite directions generating
a stable limit cycle from the focus $A.$
    \par
Fix these parameters under $\alpha_2=\alpha_2^*,$ $\beta=\beta^*$
and input the parameter $\alpha_0>0$ into (3.7) getting again
system (3.5). Since, by our assumption, this system has
two limit cycles for $\alpha_0<\alpha_0^*,$ there exists
some value of the parameter, $\alpha_0^{12}$
$(0<\alpha_0^{12}<\alpha_0^*),$ for which a semi-stable limit cycle,
$\Gamma_{12}^O,$ appears in system (3.5) and then splits into a
stable cycle, $\Gamma_1^O,$ and an unstable cycle, $\Gamma_2^O,$
on further increasing $\alpha_0.$ The formed domain $D_2^O,$
bounded by the limit cycles $\Gamma_1^O,$ $\Gamma_2^O$ and
filled by the spirals, will enlarge since, by the properties of
a field rotation parameter, the interior unstable limit cycle
$\Gamma_2^O$ will contract and the exterior stable limit cycle
$\Gamma_1^O$ will expand on increasing $\alpha_0.$ The distance
between the spirals of the domain $D_2^O$ will naturally increase,
what will prevent the appearance of a semi-stable limit cycle in
this domain for $\alpha_0>\alpha_0^{12}.$ Thus, there are no such
values of the parameters $\alpha_0^*>0,$ $\alpha_2^*<0$ and
$\beta^*>0,$ for which system (3.5) would have an additional
semi-stable limit cycle.
    \par
Obviously, there are no other values of the parameters $\alpha_0,$
$\alpha_2$ and $\beta,$ for which system (3.5) would have more
than two limit cycles surrounding the point $O$ and simultaneously
more than one limit cycle surrounding the point $A$ (on the
same reasons). It follows that system (3.5) can have at most
three limit cycles and only in the $(2\!:\!1)$-distribution.
    \par
Suppose that system (3.5) has two limit cycles, $\Gamma_1^O$ and
$\Gamma_2^O,$ around the origin $O$ and the only limit cycle,
$\Gamma_1^A,$ around the point $A.$ Fix the parameters
$\alpha_0>0,$ $\alpha_2<0,$ $\beta>0$ and input the fourth
parameter, $\gamma>0,$ into (3.5) getting system (3.1).
On increasing $\gamma,$ the vector field of~(3.1) is rotated in
positive direction, the focus $O$ changes the character of its
stability, when $\gamma=\beta-\alpha_0,$ and a stable limit cycle,
$\Gamma_3^O,$ appears from the origin, since the parameter $\alpha_2$
is non-rough and negative when $\gamma=\beta-\alpha_0$
(the Andronov--Hopf bifurcation). On further increasing $\gamma,$
the cycle $\Gamma_3^O$ will join with $\Gamma_2^O$ forming
a semi-stable limit cycle, $\Gamma_{23}^O,$ which will disappear
in a ``trajectory concentration'' surrounding the origin $O;$
the other cycles, $\Gamma_1^O$ and $\Gamma_1^A,$ will expand
disappearing in a separatrix eight-loop of the saddle $S.$
    \par
Let system (3.1) have four limit cycles: $\Gamma_1^O,$
$\Gamma_2^O,$ $\Gamma_3^O$ and $\Gamma_1^A.$ Can an additional
semi-stable limit cycle appear around the origin on increasing
the parameter $\gamma\,?$ It is clear that such a limit cycle
cannot appear either in the domain $D_2^O$ bounded by $\Gamma_1^O$
and $\Gamma_2^O$ or in the domain $D_4^O$ bounded by the origin and
$\Gamma_3^O$ because of the increasing distance between the spiral coils
filling these domains on increasing $\gamma.$ Consider two other domains:
$D_1^O$ bounded on the inside by the cycle $\Gamma_1^O$ and
$D_3^O$ bounded by the cycles $\Gamma_2^O$ and $\Gamma_3^O.$
As before, we will prove impossibility of the appearance of
a semi-stable limit cycle in these domains by contradiction.
    \par
Suppose that for some set of values of the parameters,
$\alpha_0^*>0,$ $\alpha_2^*<0,$ $\beta^*>0$ and $\gamma^*>0,$
such a semi-stable cycle exists. Return to system (3.2) again
and input first the parameters $\alpha_0>0,$ $\gamma>0$
and then the parameter $\alpha_2<0\!:$
    \vspace{-2mm}
	  $$
    \begin{array}{c}
    \dot{x}=y,\\[1mm]
    \dot{y}=-(1/2)x(x-1)(x-2)\!+\!(\alpha_0+\gamma+\alpha_2\,x^2)\,y\!
    +\,(c+dx)\,y^2+\gamma\,y^3\!.
    \end{array}	
	  \eqno(3.8)
	  $$
Fix the parameters $\alpha_0,$ $\gamma$ under the values $\alpha_0^*,$
$\gamma^*,$ respectively. On decreasing the parameter $\alpha_2,$ 
a~big stable limit cycle $\Gamma_1^{bc}$ appears from infinity and
then it contracts forming a~separatrix eight-loop of the saddle $S$
around the points $O$~and~$A.$ On further decreasing $\alpha_2,$
two stable limit cycle, $\Gamma_1^{O}$ and $\Gamma_1^{A},$ will
appear from the eight-loop surrounding $O$~and~$A,$ respectively.
Fix $\alpha_2$ under the value $\alpha_2^*$ and input the parameter
$\beta>0$ into (3.8) getting system (3.1).
    \par
Since, by our assumption, system (3.1) has three limit cycles
around the origin $O$ for $\beta<\beta^*,$ there exists
some value of the parameter, $\beta_{23}$
$(0<\beta_{23}<\beta^*),$ for which a semi-stable limit cycle,
$\Gamma_{23}^O,$ appears in this system and then it splits
into an unstable cycle, $\Gamma_2^O,$ and a stable cycle,
$\Gamma_3^O,$ on further increasing $\beta.$ The formed
domain $D_3^O$ bounded by the limit cycles $\Gamma_2^O,$
$\Gamma_3^O$ and also the domain $D_1^O$ bounded on the inside
by the limit cycle $\Gamma_1^O$ will enlarge and the spirals
filling these domains will untwist excluding a possibility of the
appearance of a semi-stable limit cycle there, i.\,e., at most
three limit cycles can exist around the origin $O.$
On the same reasons, a semi-stable limit cannot appear
around the point $A$ on increasing the parameter $\beta,$
i.\,e., at most one limit cycle can exist around this point
simultaneously with at most three limit cycles surrounding
the origin.
    \par
All other combinations of the parameters $\alpha_0,$ $\alpha_2,$
$\beta$ and $\gamma$ are considered in a similar way. It follows
that system (3.1) can have at most four limit cycles and only
in the $(3\!:\!1)$-distribution. The same conclusion can be made
for system (1.1). The theorem is proved.\quad$\Box$

\section{Application of the Wintner--Perko termination principle}

For the global analysis of limit cycle bifurcations in \cite{G},
we used the Wintner--Perko termination principle which connects
the main bifurcations of limit cycles~\cite{P}. Let us
formulate this principle for the polynomial system
    $$
    \mbox{\boldmath$\dot{x}$}=\mbox{\boldmath$f$}
    (\mbox{\boldmath$x$},\mbox{\boldmath$\mu$)},
    \eqno(4.1)
    $$
where $\mbox{\boldmath$x$}\in\textbf{R}^2;$ \
$\mbox{\boldmath$\mu$}\in\textbf{R}^n;$ \
$\mbox{\boldmath$f$}\in\textbf{R}^2$ \
$(\mbox{\boldmath$f$}$
is a polynomial vector function).
    \par
    \medskip
\noindent\textbf{Theorem 4.1
    (Wintner--Perko termination principle).}
    \emph{Any one-para\-me\-ter fa\-mi\-ly of multiplicity-$m$
limit cycles of relatively prime polynomial system \ $(4.1)$
can be extended in a unique way to a maximal one-parameter
family of multiplicity-$m$ limit cycles of \ $(4.1)$ which
is either open or cyclic.}
    \par
\emph{If it is open, then it terminates either as the parameter
or the limit cycles become unbounded; or, the family terminates
either at a singular point of \ $(4.1),$ which is typically
a fine focus of multiplicity~$m,$ or on a $($compound$\,)$
separatrix cycle of \ $(4.1),$ which is also typically of
multiplicity~$m.$}
    \medskip
    \par
The proof of the Wintner--Perko termination principle
for general polynomial system (4.1) with a vector parameter
$\mbox{\boldmath$\mu$}\in\textbf{R}^n$ parallels the proof
of the pla\-nar termination principle for the system
    $$
    \vspace{1mm}
    \dot{x}=P(x,y,\lambda),
        \quad
    \dot{y}=Q(x,y,\lambda)
    \eqno(4.2)
    \vspace{2mm}
    $$
with a single parameter $\lambda\in\textbf{R};$ see \cite{G,P}.
In particular, if $\lambda$ is a field rotation parameter of (4.1),
it is valid the following Perko's theorem on monotonic families
of limit cycles.
    \par
    \medskip
\noindent\textbf{Theorem 4.2.}
    \emph{If $L_{0}$ is a nonsingular multiple limit cycle of
$(4.2)$ for $\lambda=\lambda\,_0,$ then  $L_{0}$ belongs to a
one-parameter family of limit cycles of $(4.2);$ furthermore\/$:$}
    \par
1)~\emph{if the multiplicity of $L_{0}$ is odd, then the family
either expands or contracts mo\-no\-to\-ni\-cal\-ly as $\lambda$
increases through $\lambda\,_0;$}
    \par
2)~\emph{if the multiplicity of $L_{0}$ is even, then $L_{0}$
bi\-fur\-cates into a stable and an unstable limit cycle as
$\lambda$ varies from $\lambda\,_0$ in one sense and $L_{0}$
dis\-ap\-pears as $\lambda$ varies from $\lambda\,_0$ in the
opposite sense; i.\,e., there is a fold bifurcation at
$\lambda\,_0.$}
    \par
    \medskip
Using Theorems~4.1 and~4.2, we can give an alternative proof
of Theorems~3.1 for system (1.1), namely, we will prove
the following theorem.
    \par
    \medskip
\noindent\textbf{Theorem 4.3.}
    \emph{There exists no system $(1.1)$ having a swallow-tail
bifurcation surface of multiplicity-four limit cycles in its
pa\-ra\-meter space. In other words, system $(1.1)$ cannot have
either a multi\-plicity-four limit cycle or four limit cycles
around a singular point, and the maximum multi\-plicity or
the maximum number of limit cycles surrounding a singular
point is equal to three. Moreover, system $(1.1)$ can have
at most four limit cycles with their only possible
$(3\!:\!1)$-distribution.}
    \medskip
    \par
\noindent\textbf{Proof.} \ The proof of this theorem is carried
out by contradiction. Consider canonical systems (3.1) with three
field rotation para\-me\-ters $\alpha_0,$ $\alpha_2,$ $\gamma$ and
a semi-rotation parameter $\beta$ which is also a field rotation one
in the half-plane $x<1.$ Suppose this system has four limit cycles
around the origin $O.$ Then we get into some domain bounded by three
fold bifurcation surfaces forming a swallow-tail bifurcation surface
of multiplicity-four limit cycles in the space of the field rotation
pa\-ra\-me\-ters $\alpha_0,$ $\alpha_2,$ $\gamma$ and $\beta.$
    \par
The cor\-res\-pon\-ding maximal one-parameter family of
multiplicity-four limit cycles cannot be cyclic, otherwise there
will be at least one point cor\-res\-pon\-ding to the limit cycle
of multi\-pli\-ci\-ty five (or even higher) in the parameter space.
Extending the bifurcation curve of multi\-pli\-ci\-ty-five limit cycles
through this point and parameterizing the corresponding maximal
one-parameter family of multi\-pli\-ci\-ty-five limit cycles by
a field-rotation para\-me\-ter, according to Theorem~4.2, we will
obtain a monotonic curve which, by the Wintner--Perko termination
principle (Theorem~4.1), terminates either at the origin or on some
separatrix cycle surrounding the origin. Since we know at least the
cyclicity of the singular point \cite{zzhtt} which is equal to three,
we have got a contradiction with the termination principle stating
that the multiplicity of limit cycles cannot be higher than the
multi\-pli\-ci\-ty (cyclicity) of the singular point in which
they terminate.
    \par
If the maximal one-parameter family of multiplicity-four limit
cycles is not cyclic, on the same principle (Theorem~4.2), this
again contradicts to the result of \cite{zzhtt} not admitting
the multiplicity of limit cycles higher than three. It follows
that the maximum multi\-plicity or the maximum number of
limit cycles surrounding the origin is equal to three.
    \par
Consider other logical possibilities. For example, suppose that
system (3.1) has for $\alpha_0>0,$ $\alpha_2<0$ and $\beta>0$
three limit cycles in the $(2:1)$-distribution: two cycles around
the point $O$ and the only one around $A.$ Let us show impossibility
of obtaining additional limit cycles around the point $A$ by means
of the parameter $\gamma.$ We can suppose that a semi-stable cycle
appears around $A$ on increasing this parameter for $\gamma>0.$
Then, applying the Wintner--Perko termination principle (Theorem~4.1),
we can show that the corresponding maximal one-parameter family of
multiplicity-three limit cycles parameterized by another field rotation
parameter, e.\,g., $\alpha_2,$ cannot terminate in the focus $A,$
since it will be a rough one for $\gamma>0.$ The only additional
limit cycle in system (3.1) can appear from the focus $O$ for
the set of $\alpha_0>0,$ $\alpha_2<0,$ $\beta>0$ and $\gamma>0,$
when $\gamma=\beta-\alpha_0.$ All other possibilities, concerning also
big limit cycles from infinity, can be considered in a similar way.
    \par
Thus, we have proved Theorem~4.3 for system (1.1) giving one more proof
of Theorem~3.1 on at most four limit cycles with their only possible
$(3\!:\!1)$-distribution in this system.\quad$\Box$

\label{lastpage}
\end{document}